\documentclass[letter,12pt]{amsart}
\usepackage[numbers]{natbib}
\newcommand{\mytitle}{An introduction to univalent foundations for mathematicians}
\newcommand{\myauthor}{Daniel R. Grayson}
\newtheorem{theorem}{Theorem}[section]
 \newtheorem{lemma}[theorem]{Lemma}
 \newtheorem{axiom}[theorem]{Axiom}
\theoremstyle{definition}
    \newtheorem{definition}[theorem]{Definition}

\theoremstyle{remark}
    
\numberwithin{equation}{section}

\newcommand{\nonempty}[1]{{\Vert #1 \Vert}}
\newcommand{\bool}{\mathop{{\it Bool}}}
\newcommand{\yes}{\mathop{{\it yes}}}
\newcommand{\no}{\mathop{{\it no}}}
\newcommand{\refl}{\mathop{{\it refl}}}
\newcommand{\true}{\mathop{{\it True}}}
\newcommand{\false}{\mathop{{\it False}}}
\newcommand{\triv}{\mathop{{\it triv}}}
\newcommand{\symm}{\mathop{{\it symm}}}
\newcommand{\trans}{\mathop{{\it trans}}}
\newcommand{\weq}{\simeq}
\newcommand{\QQ}{\mathbb{Q}}
\newcommand{\ZZ}{\mathbb{Z}}
\newcommand{\NN}{\mathbb{N}}
\newcommand{\isom}{\cong}
\usepackage{mathtools}          
\newcommand{\defeq}{\vcentcolon=}

\newcommand\blfootnote[1]{%
  \begingroup
  \renewcommand\thefootnote{}\footnote{#1}%
  \addtocounter{footnote}{-1}%
  \endgroup}

\usepackage{hyperref} 
\hypersetup{pdftitle={\mytitle}}
\hypersetup{pdfauthor={\myauthor}}

\begin{document}
\title{\mytitle}
\author{\myauthor}
\address{2409 S.{} Vine St., Urbana, IL 61801, USA}
\email{drg@illinois.edu}
\urladdr{\href{http://dangrayson.com/}{http://dangrayson.com/}}
\date{Nov 4, 2017}

\begin{abstract}
  We offer an introduction for mathematicians to the univalent foundations of Vladimir Voevodsky, aiming to explain how he chose to encode
  mathematics in type theory and how the encoding reveals a potentially viable foundation for all of modern mathematics that can serve as an
  alternative to set theory.
\end{abstract}

\maketitle

\numberwithin{theorem}{section}

\tableofcontents

\section*{Introduction}

The\blfootnote{MR Author ID: 76410} traditional foundation for mathematics, chosen more than a century ago, pins down basic issues, such as what
numbers ``really are'', by giving them a specific arbitrary internal structure based on sets\footnote{For example, one may define $3$ to
  be either $\{0,1,2\}$ or $\{2\}$.}, which is irrelevant to modern daily mathematical
discourse.  Bertrand Russell's contemporaneously proposed alternative foundation \citep{MR1506041} was based on the theory of ``types'', in
which each variable is to be accompanied by a type, which is drawn from a hierarchy of types of increasing complexity and which provides the
variable's range of values, the aim being to prevent the formulation of paradoxical concepts, such as the set of all sets.  Over the intervening
decades type theory has been developed by computer scientists into a useful tool for verifying correctness and security of computer languages and by
mathematicians into a useful tool for computer verification of mathematical theorems, such the Four Color Theorem (the original 1976
computer-assisted \citep{10.2307/986491} Appel-Haken proof \citep{MR0543792,MR0543793} had gaps, but a proof was finally verified
\citep{MR2463991} in 2005 by computer using the proof assistant {\em Coq} \citep{Coq}), the Feit-Thompson odd-order theorem
\citep{10.1007/978-3-642-39634-2_14}, and the Kepler Conjecture about sphere packing \citep{kepler}.

In this brief introduction, we'll take a glance at the world of mathematics as viewed through the {\em univalent foundations} of Voevodsky,
which is based on type theory, is part of the world of {\em homotopy type theory}, is formally and precisely specified, and has been under
development over the last several years.  The most fundamental novelty for mathematicians is that the notion of equation is recast using
(dependent) types, as originally used by de Bruijn \citep{MR0274219} in the proof checker {\em Automath}, and as formulated and augmented by Per
Martin-L\"of \citep{MLTT79}, so that an equation is no longer necessarily a proposition.  Building on that, harking back to a conjecture of
Grothendieck about $\infty$-groupoids, Voevodsky singles out the types with {\em $h$-level} (at most) $n$, for a natural number $n$, with the ``true
propositions'' being the types at $h$-level $0$, with the ``propositions'' being the types at $h$-level (at most) $1$, with the ``sets'' being the types at
$h$-level (at most) $2$, with the objects of a category being the elements of a type at $h$-level 3 that captures the groupoid of isomorphisms of the category,
and so on.  (The quotation marks in the previous sentence are intended to indicate initially that these ``propositions'' and ``sets'' are
intended to play the role that propositions and sets play in set theory, but are not the same: rather, they are alternative ways to formalize a
common intuitive understanding of the terms\footnote{Perhaps one would prefer to use three terms, when needed to avoid ambiguity: {\em set} for
  the intuitive notion, {\em Zermelo-Fraenkel structure} for the set theory notion (as Voevodsky preferred to call it), and {\em $h$-set} for
  the univalent foundations notion.}.)  Thus propositions and their proofs become part of the language of mathematics, rather than part of the
language of logic\footnote{A simple way to embed logic into set theory is to associate to a proposition $P$ the set $\{x \in \{0\} \mid
  P\}$.  Under this embedding the true propositions are associated with one-element sets, propositions are associated with sets with at most
  one element, and the logical operations may be re-implemented as set theoretic operations.  Something similar happens in univalent foundations.},
and sets become a special case of something more fundamental, {\em types}, worthy of independent study.

The formal mathematical language, together with the Univalence Axiom, fulfills the mathematicians' dream: a language for mathematics invariant under
``equivalence'' and thus freed from irrelevant details and able to merge the results of mathematicians taking different but equivalent
approaches.  Voevodsky called this invariance property of the language ``univalence''\footnote{Voevodsky explained his choice of the term
  ``univalent'' in \citep{VV-IHP-talk}: it comes from a Russian translation of the Boardman and Vogt book ``Homotopy Invariant Algebraic
  Structures on Topological Spaces'', where the term ``faithful functor'' is translated as ``univalent functor''.  He also said ``Indeed these
  foundations seem to be faithful to the way in which I think about mathematical objects in my head.''}.
It offers the hope that formalization and verification
of today's mathematical knowledge may be achievable, relieving referees of mathematics articles of the tedious chore of checking the details of
proofs for correctness, allowing them to focus on importance, originality, and clarity of exposition.  Feasibility and practicality of the
approach were demonstrated by Voevodsky in his {\em Foundations} \citep{Foundations,UniMath2015}, based on {\em Coq}, and various teams have continued
formalization efforts based on it.  In section \ref{form} we describe how the notions of {\em group} and {\em torsor} are encoded in the system,
thereby motivating the definitions of {\em proposition} and {\em set} and aiming to expose the aspects of the system that make the formalization
succeed.

Further useful expositions of the subject include these: \citep{MR3363596,hottbook,Swansea-Shulman}, and some philosophical background is provided in
\citep{Tsem1,Tsem2,Tsem3,Tsem4}.

I turn now to a bit of speculation.  The beauty of mathematics, as appreciated by human mathematicians, is enhanced by formalization, for a
beautiful proof known to be correct is more beautiful than the same proof in unverified form.  Even if formalization becomes widely used, it
will still be humans who decide what to prove, how to prove it, which proofs are better than others, and which proofs are worth publishing in
journals.  A beautiful incomplete or incorrect proof can benefit from formalization, too -- then anyone will know precisely what has been done,
that it has been done correctly, and what remains to be done.

The adoption of proof formalization by a broad segment of practicing mathematicians will await the day when mathematicians feel more successful
with formalization than without, as happened with {\em TeX} starting in the 1980's.  There is much development to be done to make proof
assistants more powerful and easier to use.  Formalization of a result also depends on formalization of the results cited by it, so another
prerequisite for wide adoption is the formalization of a large body of existing mathematics.

Current proof assistants are oriented toward enabling humans to enter proofs efficiently, so perceiving the beauty of a proof by reading its
formalization will not be as enlightening as by reading human-written prose about it, except for those who become proficient with the
proof-entering mechanism.  Eventually, when proof-entering technology is sufficiently developed, we can expect proof-reading technology to be
developed further, and, ultimately, the optimal way to read a proof and to learn it will be with the active assistance of a proof assistant, so
that details can be clarified upon demand, right down to the bottom.  Anyone who has worked their way through even a fragment of EGA, with its
numerous cross-references to earlier results, knows how desirable that is in a text, and how hard it is to achieve outside a unified framework.
In an automated comprehensive system, it would be effortless.

Alexander Beilinson has written \cite{IMGelfand}, perhaps expressing a thought of I.{} M.{} Gelfand: ``Modern mathematics is a unique thrust of
conceptual thought: once the right concept (a mathematical structure) and a language to deal with it are found, a whole new world unfolds.''
This is what has happened with univalent foundations, and in this new world mathematicians and the beauty of mathematics will prosper.

\section*{Dedication}

I dedicate this article to Vladimir Voevodsky, an ingenious mathematician with an immense drive to create, who died suddenly in September, 2017,
at the age of 51.
Our friendship began perhaps in 1994, when I first became aware of his work on motivic cohomology, and he visited me in Urbana to explain what
it was about.  Some mathematicians had suspected that there was a larger role for homotopy theory to play in motivic cohomology, but it was up
to Vladimir to show the way, which he did with immense energy, and a detailed plan for the future covering many components of his program.  His
work, including a proof of the Milnor Conjecture (1970), won him the Fields Medal in 2002, went most of the way toward the ultimate solution of
the Bloch-Kato Conjecture (of which the Milnor Conjecture is a special case), the Beilinson-Lichtenbaum Conjecture, and the conjecture of Milnor
about the Witt ring of quadratic forms.  The field of motivic homotopy theory remains a vibrant one to this day, although Vladimir stopped
participating in it about ten years ago or so.  He started thinking about the use of computers for formal representation and verification of
mathematical statements and proof in 2002.  We worked together for a week in Spring, 2004, brainstorming about that, dreaming about what an
ideal system would be like, and I even wrote some code.  I soon became occupied with other matters, but he surveyed the field, and by 2006 had
chosen type theory as the proper language.  By 2010 he had chosen a suitable encoding of mathematics in type theory (described in this article)
and had spent three months writing some brilliant, beautiful, and instructive proofs in {\em Coq}, which he dubbed {\em Foundations} \citep{Foundations}, later
incorporated into a broader project with more authors called {\em UniMath} \citep{UniMath}.  Always eager to participate in his project, I finally found the time to start collaborating with him in 2011,
and in 2012-2013 I attended the special year on the topic at the Institute for Advanced Study and found it bracing and enlightening.  His final
project was to establish the soundness of the univalent foundations with complete mathematical rigor in a series of papers, eight of which he managed to write (see section
\ref{the_interpretation}).  His dream to establish the univalent foundations as the preferred and practical framework for formalizing the
world's mathematical knowledge seems feasible, but much work remains to be done by the community.  A fitting memorial for Vladimir would be to
formalize his work on motivic homotopy theory in the univalent foundations.

\section*{Acknowledgments}

I thank the {\em Oswald Veblen Fund} and the {\em Friends of the Institute for Advanced Study} for supporting my stay at the Institute for
Advanced Study in Winter/Spring, 2017, where I worked on part of this paper, and the Institute for Advanced Study School of Mathematics {\em
  Fund for Mathematics} for supporting my stay in Winter/Spring, 2018, where I finished it.  I thank Benedikt Ahrens, Guillaume Brunerie, 
Thierry Coquand, Bas Spitters,  Mart\'in Escard\'o, and Anders M\"ortberg for useful comments on earlier drafts of this paper.

\section{What is a type?}

In some computer programming languages, all variables are introduced along with a declaration of the type of thing they will refer to.  For example,
one may encounter types such as $bool$, $string$, $int$, and $real$, describing Boolean values, character strings, 32 bit integers, and 64 bit
floating point numbers.  The types are used to determine which statements of the programming language are grammatically well-formed.
For example, if $s$ of type $string$ and $x$ is of type $real$, we may write $1/x$, but we may not write $1/s$.

Types occur in traditional mathematics, as expressed in informal mathematical speech, and are used in the same way: all variables are introduced
along with a declaration of the type of thing they will refer to.  For example, one may say ``consider a group $G$'', ``consider a ring $R$'',
or ``consider an algebraic variety $X$ over a field $k$''.  One does not see circumlocutions such as ``consider a thing $X$: if $X$ is a group
then \dots''.

This informal use of types is not supported by first-order logic and Zermelo-Fraenkel set theory, where there are just two types of mathematical
or logical objects: propositions and sets.  Nevertheless, as with computer languages, here the types are used to determine which statements of
the theory are grammatically well-formed.  For example, if $X$ and $Y$ are sets, then we may write the proposition $X \in Y$, and if $P$ and $Q$
are propositions we may write the proposition $P \wedge Q$, but we may not write $P \in Q$ or $X \wedge Y$.

In type theory, there are many more types, they support the traditional use of types in informal mathematical speech, and they are used for
everything.  There will be enough primitive ways to form new types from old ones to provide everything we need to formalize mathematics,
including: logical quantifiers, functions, families, pairs, products, sums, and equality.  Building on that, one may introduce types such as
$\NN$, $\ZZ$, $\QQ$, $Group$, $Ring$, and $AbelianCategory$, whose elements are, respectively, the natural numbers, the integers, the rational
numbers, the groups, the rings, and the abelian categories.

One expresses the statement that an ``element'' $a$ is of ``type'' $X$ by writing $a:X$.  Using that notation, each variable is introduced along
with a declaration of the type of thing it will refer to, and the declared types of the variables are used to determine which statements of the
theory are grammatically well-formed.  Propositions and sets are to be re-implemented as types with a certain property, and we will have no use
for the set membership operator $\in$.

As we have said, if $X$ and $Y$ are types, there will be a type whose elements serve as {\em functions} from $X$ to $Y$; the notation for it is
$X \to Y$.  This allows us to introduce the surprising mathematical pun $f : X \to Y$, which says that $f$ is an element of the type $X \to Y$,
and which can be read traditionally as saying that $f$ is a function from $X$ to $Y$.\footnote{This new interpretation of the notation may
  appear jarring to those who are accustomed to writing something like ``consider a function $X \xrightarrow f Y$'', regarding the arrow as a
  pictorial representation of the function itself; we won't do that.}

Functions behave as one would expect, and one can make new ones in the usual way.\footnote{ To provide an example of making new functions in the
  usual way, consider functions $f : X \to Y$ and $g : Y \to Z$.  We define their composite $g \circ f : X \to Z$ by setting $g \circ f \defeq
  (a \mapsto g(f(a)))$.  Such definitions are to be regarded as syntactically transparent in our formal system, in the sense that two formal
  expressions will be regarded as being {\em the same by definition} if they yield the same formal expression after the definitions of all the
  symbols within them are completely expanded.  Given two expressions that are the same by definition, we may replace one with the other in any
  other expression, at will.  Here is an example: consider functions $f : X \to Y$, $g : Y \to Z$, and $h : Z \to W$.  Then $(h \circ g) \circ
  f$ and $h \circ (g \circ f)$ are the same by definition, since applying the definitions within expands both to $a \mapsto h(g(f(a)))$.

  One may define the identity function $id_X : X \to X$ by setting $id_X \defeq (a \mapsto a)$.  Application of definitions shows that $f \circ
  id_X$ is the same as $a \mapsto f(a)$, which, by a standard convention, is to be regarded as the same as $f$.  A similar computation applies
  to $id_Y \circ f$.
}

In the following sections we will expose various other elementary sorts of types, focusing most on those with novel aspects.

\section{The type of natural numbers}\label{nat}

The inductive definition of the natural numbers provided by Giuseppe Peano in 1889 \citep{peano-principia} is mirrored almost exactly in type
theory.

Here are Peano's rules\footnote{There is a reason for not calling them Peano's Axioms in this context, since we wish to reserve the word ``axiom''
  for decrees that a proposition has a proof.} for constructing the natural numbers in the form that is used in type theory.
\begin{enumerate}
\item[P1:] there is a type called $\NN$ (whose elements will be called natural numbers);
\item[P2:] there is an element of $\NN$ called $0$;
\item[P3:] if $m$ is a natural number, then there is also a natural number $S(m)$, called the {\em successor} of $m$;
\item[P4:] given a family of types $X(m)$ depending on a parameter
  $m$ of type $\NN$, in order to define a family $f(m) : X(m)$ of elements of each of them it suffices to provide an element $a$ of $X(0)$ and
  to provide, for each $m$, a function $g_m : X(m) \to X(S(m))$.  (The resulting function $f$ may be regarded as having been defined inductively
  by the two declarations $f(0) \defeq a$ and $f(S(m)) \defeq g_m(f(m))$.)
\end{enumerate}
\nopagebreak
Mathematicians will recognize rule P4 as ``the principle of mathematical induction'' when each $X(m)$ is a proposition, and as ``defining a
function by recursion'' if each $X(m)$ is a fixed set $Y$.  We will refer to it simply as ``induction for $\NN$''.\footnote{Here is an example
  of defining a function by recursion using induction for $\NN$.  We define the factorial function $f : \NN \to \NN$ by setting $f(0) \defeq 1$
  and setting $f(S(m)) \defeq (m+1) \cdot f(m)$.  One can infer that the function $g_m$ of rule (P4) is $n \mapsto (m+1) \cdot n$.}

Notice that
the two cases in an inductive definition correspond to the two ways of introducing elements of $\NN$ via the use of rules P2 and P3.
Intuitively, the induction principle for $\NN$ amounts to saying that the element $0$ and the function $S$ ``generate'' the type $\NN$, in the
same way that a set may generate the corresponding free group.

We introduce the following definitions.
\begin{align*}
 1 & \defeq S(0) \\
 2 & \defeq S(1) \\
 3 & \defeq S(2) \\
 4 & \defeq S(3)
\end{align*}

We may use induction to define the sum $m+n$ of two natural numbers, as a natural number.  We handle the two possible cases for the argument $m$
as follows: we define $0+n \defeq n$, and we define $(S(m))+n \defeq S(m+n)$.  Application of definitions shows, for example, that $2+2$ and $4$
are the same by definition, because they both reduce to $S(S(S(S(0))))$.\footnote{The reduction of the expression $2+2$ to $S(S(S(S(0))))$ by
  applying the definitions involved can be regarded as a computation, which is ``complete''.  Computation succeeds in other cases involving
  natural numbers, too.  That is why we don't refer to Peano's rules as axioms, as we mentioned in a previous footnote, for we wish to reserve
  the word ``axiom'' for statements unaccompanied by effective methods for computation.  A good example of such an axiom is the law of the
  excluded middle, which states, in our context, that any proposition $P$ either has a proof or its negation has a proof.  The law provides no
  effective way to decide which case we are in nor what the proof would be in that case.  Nevertheless, it is consistent with type theory and
  all the axioms we intend to use, so its use as a hypothesis in a theorem is permissible.  If one wishes to define a discontinuous function of
  a real variable, it is needed, because there is no effective way to decide whether a real number $x$ is, for example, positive -- even a very
  long computation may yield a very accurate rational approximation to $x$ that is too close to $0$ to determine whether $x$ is positive.}

Before we can write an equation such as $2+2=4$, we must introduce a formal treatment of equality in type theory.  We do that in the next section.

\section{Identity types}\label{paths}

The most important type is the {\em identity type}, which implements the intuitive notion of equality, but in a novel way; the mathematical
reader will be more comfortable if we call it the {\em equality type}, at least initially.  Its novelty derives from three aspects: equality
between two elements may be considered only when the two elements are of the same type; equality may fail to be a proposition (unless it's an
equality between two elements of a set); and a proof of equality between two sets will amount to a bijection between them.

We can anticipate the consequences of the first of those three aspects immediately.  For example, let $2_\ZZ$ denote the integer corresponding
to $2$ and let $2_\QQ$ denote the rational number corresponding to $2$.  Then there is no way to translate the valid but mathematically
irrelevant equation $2_\ZZ = 2_\QQ$ of set theory into type theory, because $2_\ZZ$ and $2_\QQ$ are elements of distinct types: $\ZZ$ and $\QQ$.
Beginning mathematics students would be completely comfortable with this restriction.  For another example, consider two sets (or types) $X$ and
$Y$.  There is no way to state a condition that says that $X$ is a subset (or subtype) of $Y$.  Instead, as in category theory, one considers
embeddings of $X$ into $Y$.

Observe that the smallest reflexive relation on the elements of a set $X$ is equality.  Per Martin-L\"of realized in the 1970's that one can use
this observation to provide an inductive definition of equality\footnote{See \citep[section 1.7]{MR0387009}.  An antecedent was published in
  \citep[section 3.8.2, p.{} 190]{MR0387023}.} analogous to the inductive definition of the natural numbers presented in section \ref{nat}.

Here are Martin-L\"of's rules for constructing equality types.
\begin{enumerate}
\item[E1:]
  for any type $X$ and for any elements $a$ and $b$ of it, there is a type $a=b$;\footnote{
    We point out the order of quantifiers in rule E1 is grammatically required.  Suppose one tried to rephrase rule E1 to have the form
    ``for any elements $a$ and $b$, \dots''.  Then the puzzle would be how to formulate the condition that $a$ and $b$ have the same type, and there
    is no way to do that in our formal language, as currently designed.  An arbitrary element of an arbitrary type is always introduced by %
    introducing a quantifier first for its type, so the quantifier for the element can state the type of the element, as required.
  }
\item[E2:] for any type $X$ and for any element $a$ of it, there is an element $\refl(a)$ of type $a=a$ (the name $\refl$ comes from the word
  ``reflexivity'')
\item[E3:] for any type $X$ and for any element $a$ of it, given a family of types $P(b,e)$ depending on parameters $b$ of type $X$ and $e$ of type
  $a=b$, in order to define elements $f(b,e) : P(b,e)$ of all of them it suffices to provide an element $p$ of $P(a,\refl(a))$.  The resulting
  function $f$ may be regarded as having been completely defined by the single definition $f(a,\refl(a)) \defeq p$.
\end{enumerate}

An element of $a=b$ can be thought of, for now, as a proof that $a$ is equal to $b$.

We see from rule E2 that $\refl(S(S(S(S(0)))))$ serves as a proof of $2+2=4$, as do $\refl(4)$ and $\refl(2+2)$.  A beginning student might wish
for a more detailed proof of that equation, but as a result of our convention above that definitions are syntactically transparent, the
application of definitions, including inductive definitions, is regarded as a trivial operation.

We will refer to rule E3 as ``induction for equality''.  It says that to prove something about (or to construct something from) every proof that
$a$ is equal to something else, it suffices to consider the special case where the proof is the trivial proof that $a$ is equal to itself, i.e.,
where the proof is $\refl(a) : a=a$.  Notice that the single case in such an induction corresponds to the single way of introducing elements of
equality types via rule E2, and compare that with P4, which dealt with the two ways of introducing elements of $\NN$.
Intuitively, the induction principle for equality amounts to saying that the element $\refl(a)$ ``generates'' the system of types $a=b$, as $b$
ranges over elements of $A$.

The mathematical reader who accepts Peano's rules as valid may feel uneasy about the validity of Martin-L\"of's rules for equality types.
Peano's rules are easily accepted as valid, because they describe things we know about, natural numbers, and they posit inference rules which
appear to be valid.  Proofs of equality, on the other hand, are not things we know about, at least, not as mathematical objects in their own
right.  The informal way to justify them is to repeat that the smallest reflexive relation on the elements of a set $X$ is equality, so
everything true about equality ought to flow from reflexivity.  Formally, one way to proceed is to justify the consistency of the formal system
by providing an interpretation of the entire formal system in a suitable mathematical structure.  The justification of type theory with the
univalence axiom in that way one of the goals of a project of Voevodsky; we cite references in section \ref{the_interpretation}.

For now, until we encounter univalence and its consequences, the reader may consider the following hints pointing toward an interpretation (in
classical mathematics) of this formal system.  In that interpretation, every type $X$ is interpreted as a set $|X|$, and an element $a:X$ is interpreted as
an element $|a| \in |X|$.  A function type $X \to Y$ is interpreted as the set of functions from $|X|$ to $|Y|$.  The type $\NN$ is interpreted
as the set $|\NN|$ of natural numbers.  The type $a=b$ is interpreted as the one-point set $\{*\}$ if $|a| = |b|$, and is interpreted as the
empty set $\{\}$ if $|a| \ne |b|$.  The element $\refl(a) : a=a$ is interpreted as the element $*$ of the one-point set.  Using this
visualization may provide some temporary intuition for Martin-L\"of's rules.  The interpretation is not required to be faithful, except in the
sense that the empty type (to be introduced below) should be interpreted as the empty set.  Then if the formal system allowed an absurdity to be
proved, the corresponding element $a$ of the empty type would give an element $|a|$ of the empty set, providing a contradiction in set theory.

Later, when univalence is introduced, the mathematical interpretation will have to be modified so types are interpreted not as sets, but as
topological spaces, to handle the behavior of types that are not sets.

We may use induction to prove symmetry of equality.  In accordance with our discussion of implication above, we show how to produce an
element of $b=a$ from an element $e$ of $a=b$, for any $b$ and $e$.  By induction (letting $P(b,e)$ be $b=a$ in rule E3 above), it suffices to produce an element of
$a=a$; we choose $\refl(a)$ to achieve that.

Transitivity of equality is established the same way.  For each $a,b,c:X$ and for each $p:a=b$ and for each $q:b=c$ we want to produce an
element of type $a=c$.  By induction on $q$ we are reduced to the case where $c$ is $b$ and $q$ is $\refl(b)$, and we are to produce an element
of $a=b$.  The element $p$ serves the purpose.  Notice the similarity of this inductive definition with the definition given above of the sum
$m+n$.

Associativity of transitivity is established the same way.  We leave its proof as an exercise.

Now let's consider how we might formulate our symmetry result.  One way would be to assert that we have proved the following lemma.

\begin{lemma}
  For any type $X$ and for any $a,b:X$, $(a=b) \to (b=a)$.
\end{lemma}

Aside from using $\to$ instead of $\implies$ for the ``implication'', it looks traditional, but there is a drawback that is perhaps not
immediately apparent: standard practice in current mathematical prose is to regard the statement of a lemma as sufficient for application of the
lemma, and to regard the proof (once it has been provided and verified) as irrelevant for the application of the lemma\footnote{I don't mean to
  imply that the proof is useless for the mathematical reader, who may learn something from it, or who may modify it to prove something else.}.
In this case, the construction of the symmetry function that we have provided above is not revealed in the statement of the lemma, but only in
the proof of the lemma.  Having the precise function available for perusal might be important if $b=a$ has more than one proof.  Thus a better
formulation would be as a definition, as follows.

\begin{definition}
  For any type $X$ and for any $a,b:X$, let $\symm_{a,b} : (a=b) \to (b=a)$ be the function defined by induction by setting 
  $\symm_{a,a}(\refl(a)) \defeq \refl(a)$.
\end{definition}

Similarly, transitivity is best formulated as an inductive definition $\trans_{a,b,c} : (a=b) \to ((b=c) \to (a=c))$.  We may abbreviate
$(\trans_{a,b,c}(p))(q)$ as $p*q$.

Given 4 points $a,b,c,d:X$ and 3 proofs $p:a=b$, $q:b=c$, $r:c=d$ of equality, we may produce two proofs that $a=d$, namely $(p*q)*r$ and
$p*(q*r)$.  Because they are elements of the same type, we may consider the type $(p*q)*r = p*(q*r)$.  An element of it is a proof of
associativity, and it may be constructed by induction.

One frequent use of equality proofs is in {\em substitution}.  Let $X$ be a type, and let $P(x)$ be a family of types depending on a parameter
$x:X$.  Suppose $a,b:X$ and $e:a=b$.  Then there is a function of type $P(a) \to P(b)$.  To prove that, it suffices, by induction, to consider
the case where $b$ is $a$ and $e$ is $\refl(a)$ and to provide a function $P(a) \to P(a)$.  To achieve that, we provide the identity function.

In the case where each $P(a)$ is a proposition (to be defined later), we may say that the truth of $P(a)$ is invariant under substitution.  In
other words, equal elements of $X$ have the same properties.  In the case where each $P(a)$ is not assumed to be a proposition, and thus its
elements are ways to add extra structure to $a$, we may say that the structures can be {\em transported} from $a$ to $b$ by the proof $e$ of
equality.

\section{Other types}

There are other examples of types that are conveniently presented as inductive definitions, in the style we have seen with the natural numbers
and the equality types.  We present three examples.

Firstly, there will be the ``empty'' type, called $\emptyset$, defined inductively, with no way to construct elements provided in the
inductive definition.  The inductive principle for $\emptyset$ says that to prove something about (or to construct something from) every element
of $\emptyset$, it suffices to consider no special cases (!).  Hence, every statement about an arbitrary element of $\emptyset$ can be proven.
As an example, we may prove that any two elements $x$ and $y$ of $\emptyset$ are equal by using induction on $x$.

An element of $\emptyset$ will be called an {\em absurdity}, and the negation $\neg P$ of a proposition $P$ will be implemented as the function
type $P \to \emptyset$.  This is sensible, because an element of $\neg P$ could be applied to an element of $P$ to produce an element of
$\emptyset$, i.e., an absurdity.

Another appropriate name for $\emptyset$ is $\false$.

We may also construct a function $\false \to X$, for any type $X$, by induction, showing that from an absurdity anything follows.

To encode the property that $X$ has no elements we use the type $X \to \emptyset$.  To encode the property that elements $a,b:X$ are not equal,
we use the type $(a=b) \to \emptyset$, and we let $a \ne b$ denote it.

Secondly, there will also be a type called $\true$, defined inductively and provided with a single element $\triv$; (the name $\triv$ comes from the word
  ``trivial'').  Its induction principle
states that, in order to prove something about (or to construct something from) every element of $\true$, it suffices to consider the special
case where the element is $\triv$.  As an example, we may prove, for any element $u : \true$, that $u=\triv$, by using induction to reduce
to proving $\triv=\triv$, a proof of which is provided by $\refl(\triv)$.  One may also prove that any two elements of $\true$ are equal by using induction twice.

There is a function $X \to \true$, for any type $X$, namely: $a \mapsto \triv$.  This corresponds, for propositions, to the statement that an
implication holds if the conclusion is true.

The name of $\true$ is appropriate, because if $P$ is a proposition, a function of type $\true \to P$ could be applied to the element $\triv$,
yielding an element of $P$, thus proving $P$.

Thirdly, there will be a type called $\bool$, defined by induction and provided with two elements, $\yes$ and $\no$.  One may prove by induction
that any element of $\bool$ is equal to $\yes$ or to $\no$.

We may use substitution to prove $\yes \ne \no$.  To do this, we introduce a family of types $P(b)$ parametrized by a variable $b:\bool$.
Define $P(\yes) \defeq \true$ and define $P(\no) \defeq \false$.  The definition of $P(b)$ is motivated by the expectation that we will be able to prove
that $P(b)$ and $b = \yes$ are equivalent.  If there were an element $e: \yes = \no$, we could substitute $\no$ for $\yes$
in $\triv : P (\yes)$ to get an element of $P(\no)$, which is absurd.  Since $e$ was arbitrary, we have defined a
function $(\yes=\no) \to \emptyset$, establishing the claim.

In the same way, we may use substitution to prove that successors of natural numbers are never equal to $0$, i.e., for any $n:\NN$ that
$0 \ne S(n)$.  To do this, we introduce a family of types $P(i)$ parametrized by a variable $i:\NN$.  
Define $P$ recursively by specifying that $P(0) \defeq \true$ and 
$P(S(m)) \defeq \false$.  The definition of $P(i)$ is motivated by the expectation that we will be able to prove that $P(i)$ and $i = 0$ are
equivalent.  If there were an element $e: 0 = S(n)$, we could substitute $S(n)$ for $0$ in $\triv : P ( 0 )$ to get an element of $P(S(n))$,
which is absurd.  Since $e$ was arbitrary, we have defined a function $(0=S(n)) \to \emptyset$, establishing the claim.

Type theory provides {\em sums} of types.  By this we mean if $X$ is a type and $Y(x)$ is a family of types indexed by a parameter $x$ of
type $X$, then there will be a type $\sum _{x:X} Y(x)$ whose elements are the pairs $(a,b)$, where $a:X$ and $b:Y(a)$.  This is reminiscent of the
sum in set theory, which is implemented by the disjoint union $\amalg _{x:X} Y(x)$.  Sums may be implemented by an inductive definition,
and they behave as one would expect.\footnote{\dots except perhaps when $X$ is a type that is not a set, because then
the type $X$ behaves like a topological space, the family $Y(x)$ behaves like a family of spaces varying continuously in the parameter $x$,
and the sum behaves like the total space of the family.}

Using sums we may encode the fibers of a function.  Given a function $f : X \to Y$ and an element $y:Y$, the fiber (or inverse image) $f^{-1}(y)$ consists of
points $x$ such that $f(x) = y$.  This is encoded by defining $f^{-1}(y) \defeq \sum_{x:X} (f(x) = y)$.  In other words, a point of the fiber is a
pair $(x,e)$ consisting of the point $x$ and a proof $e$ of the equation $f(x) = y$.

There is a binary sum operation on types: for any types $X$ and $Y$, there is a type $X \amalg Y$, which is analogous to the disjoint union of
two sets.  From an element of $X$ or an element of $Y$ we can produce an element of $X \amalg Y$.  The binary sum can be implemented as a sum
where the index type is $\bool$ and the family of types sends $\yes$ to $X$ and sends $\no$ to $Y$, or it may be implemented as an inductive
definition.  Binary sums behave as one would expect.

Our type theory will also contain {\em products} of types.  By this we mean if $X$ is a type and $Y(x)$ is a family of types indexed by a
parameter $x$ of type $X$, then there will be a type $\prod _{x:X} Y(x)$ whose elements serve as {\em families} of elements $b:Y(a)$, one for each $a:X$.
A family is much like a function, but without a single codomain.  This is reminiscent of the product in set theory, which uses the same notation.
Products behave as one would expect.\footnote{\dots except perhaps when $X$ is a type that is not a set, because then the product behaves like
  the space of continuous sections of the family.}

There is a binary product operation on types: for any types $X$ and $Y$, there is a type $X \times Y$.  From an element of $X$ and an element of
$Y$ we can produce an element $(x,y)$ of $X \times Y$.  The binary product can be implemented as an inductive definition, as a special case of
products, or as a special case of sums.  Binary products behave as one would expect.

\section{Formalization of mathematics}\label{form}

With equality types available, with all their expected properties, we may encode some elementary mathematical properties as types, to show how such
encoding goes in practice, as implemented (approximately) in: the {\em UniMath} project \citep{UniMath}, which is exposed by Voevodsky
in \citep{UniMath2015}; as in the {\em HoTT} project \citep{HoTT}, which is exposed in \citep{1610.04591}; as in the {\em HoTT-Agda} project
\citep{HoTT-Agda}; and as in the {\em Lean} theorem prover \citep{Lean}.

Let $X$ be a type.  The property that $X$ has at most one element is equivalent to the property that any two
elements are equal, so is encoded by $\prod_{a:X} \prod_{b:X} (a=b)$.  The property that $X$ has exactly one element is equivalent to having an
element such that every other element is equal to it; hence it is encoded by $\sum_{a:X} \prod_{b:X} (a=b)$.

Now let us consider how to formalize the mathematical notion of group.  More precisely, we adopt as our goal the definition of a type, $Group$,
whose elements are the groups, whose existence was advertised in the introduction.

Following a traditional approach, one may choose to define a group as a tuple $(G,e,i,m)$, where $G$ is a set, $e$ is the unit element, $i$ is
the inverse operation, and $m$ is the multiplication operation; one also asserts the truth of various properties of the operations, such as
associativity.

Using sums, we may define a 4-tuple $(G,e,i,m)$ as an iterated pair of the form $(G,(e,(i,m)))$, and we may consider defining $Group$ to be the
type of such 4-tuples.

Well, almost: we need to specify a type whose elements will provide the candidates for $G$, as the notion of sum requires.  For that purpose, we
introduce a {\em universe}, $U_0$; its elements will correspond to ``small'' types.  The notion of ``universe'' seems to have first arisen in
the formalization of set theory by Bernays \citep{bernays-1}, building on work of von Neumann, in which there are three types of thing:
propositions, sets, and classes.  One of the classes has all of the sets as its elements, and it is analogous to our $U_0$.  One could extend
their system by introducing {\em hyperclasses}, one of which would have all of the classes as its elements, and we could regard it as analogous to a
bigger universe $U_1$, one of whose elements is $U_0$.  And so on.  Indeed, such a chain of universes $U_0, U_1, U_2, \dots$ is postulated in our type theory, and the
concept turns out to be of fundamental importance.  Each one is a ``universe'' in the sense that it is closed under the operations of type
theory that make new types from old ones.\footnote{The interpretation in classical mathematics that we discussed above for verifying consistency
  of our formal system will need to be adjusted to provide interpretations of the universes.  For that purpose one assumes that there is an
  ascending chain of Grothendieck universes $|U_0| \in |U_1| \in |U_2| \in \cdots$, and one trusts that that assumption doesn't add an inconsistency to set
  theory.}

Having introduced universes, we may say that our tuples $(G,e,i,m)$ will be the elements of type $\sum _ {G : U_n} \sum _ {e:G} \sum _ {i:G
  \to G} ( G \times G \to G )$.  Here $n$ is a fixed universe level, and we are engaged in specifying a type to encode the groups that live in
$U_n$.

In order to formalize the assertion that the operations satisfy various properties, we recall that all such propositions are to be
re-implemented as types and their proofs are to be re-implemented as elements of those types.  So we extend our 4-tuple to a 9-tuple
$(G,e,i,m,\alpha,\lambda,\rho,\lambda',\rho')$, where $\alpha$ is a proof of associativity, $\lambda$ is a proof of the left unit law, $\rho$ is
a proof of the right unit law, $\lambda'$ is a proof of the left inverse law, and $\rho'$ is a proof of the right inverse law.

Well, almost.  The associativity property is formulated in set theory as
$$\forall _ {a \in G}\forall _ {b \in G}\forall _ {c \in G} m (m (a,b), c) = m(a, m(b,c)),$$ so we need a type that can represent universal quantifiers.
The product type will serve that purpose, so we can use
 $$\prod _ {a: G}\prod _ {b: G}\prod _ {c: G} m (m (a,b), c) = m(a, m(b,c)).$$ We will take $\alpha$ to be an element of that type; in other
words, it will be a function that provides, for each $a,b,c$ a proof of the equation $m (m (a,b), c) = m(a, m(b,c))$.  The other four properties
are implemented in a similar way.

Now consider the possibility that we have two unequal proofs, $\alpha$ and $\alpha'$, of associativity.  Then we would get two unequal groups,
$$(G,e,i,m,\alpha,\lambda,\rho,\lambda',\rho')$$ and $$(G,e,i,m,\alpha',\lambda,\rho,\lambda',\rho'),$$ and that would be an unintended departure from
traditional mathematics.  This is a consequence of the necessity to include the proofs of the properties in the tuple, because we are trying to
define the type whose elements are the groups in our formal language, which provides no way to assert the properties of a group ``on the side''.
Thus we are led to insist that there is a proof of $\alpha = \alpha'$, as a consequence of which it would follow that our two groups are equal.
Each of the functions $\alpha$ and $\alpha'$ assigns to any $a,b,c:G$ a proof of $m (m (a,b), c) = m(a, m(b,c))$, and functions with unequal
values can't be equal, so we will need to be able to prove that any two proofs of $m (m (a,b), c) = m(a, m(b,c))$ are equal.  For this purpose
it would be most convenient if it were true more generally that any two proofs of an equation $x=y$ between any two elements $x$ and $y$ of $G$ are
equal.

Thus we are led to embark on a diversion about types for which any two elements are equal.  In informal mathematical speech and in set theory,
any proof of a proposition is irrelevant, as it tells us only that the proposition is true: it provides no information that is needed later on.
Good mathematical prose is formulated to follow that practice by relocating information needed later on from the proof to the statement
being proven.  That and the preceding paragraph motivate our definition of ``proposition'' in this context.  A {\em proposition} is a type such
that any two elements of it are equal.  The statement that a type $P$ is a proposition is encoded by the type $isProp(P) \defeq \prod_{a:P}
\prod_{b:P} (a=b)$.\footnote{For the purpose of exposition, we have chosen a definition that differs from Voevodsky's and is simpler, but equivalent
  to it.}

Statements proven previously can be rephrased as saying that $\false$ and $\true$ are propositions.

Suppose $P$ and $Q$ are propositions.  The implication $P \implies Q$ can be implemented as the type $P \to Q$.  An element of it will transform
any proof of $P$ into a proof of $Q$.  The conjunction $P \wedge Q$ can be implemented as the product $P \times Q$.
We may define the universal quantifier $$\forall _{x:X} P(x) \defeq \prod _{x:X} P(x)$$ for a family of propositions $P(x)$ depending on a parameter
$x:X$.  It is a proposition, and it has the properties one would expect.
By contrast, the disjunction
$P \vee Q$ cannot be implemented as the binary sum $P \amalg Q$, because that type will have two elements if both $P$ and $Q$ have proofs.  We
will return to this issue later.  A similar remark applies to an arbitrary sum of propositions: it cannot be used to implement the
existential quantifier.

We continue our diversion with a discussion of types $X$ with the property that any two proofs of an equation $a=b$ between any two elements $a$
and $b$ of $X$ are equal, motivated by our need for $G$ to be such a type, as we said above.  With our definition of proposition in hand, we can rephrase the
condition on $X$ as positing that the equality type $a=b$ between any two elements $a,b:X$ is a proposition.  A type $X$ with that property is
defined by Voevodsky to be a ``set''.  Let us use the notation $isSet(X)$ to denote the type encoding this condition.

It can be proven that the propositions are the sets with at most one element, that many types, including $\bool$ and $\NN$, are sets, and that
many operations for making new types from old preserve the property of being a set.

Ending our diversion and returning to our formalization of the notion of ``group'', we see that we need $G$ to be a set.  Since we can't
postulate that $G$ is a set ``on the side'', we have to add a proof of it to the tuple.  Thus our final definition of ``group'' is that it is a
10-tuple $$(G,e,i,m,\alpha,\lambda,\rho,\lambda',\rho',\iota),$$ where $\iota$ is a proof that $G$ is a set, i.e., it is an element of the type
$isSet(G)$.

Now consider the possibility that we have two unequal proofs, $\iota$ and $\iota'$, that $G$ is a set.  Then we would get two unequal groups,
$$(G,e,i,m,\alpha,\lambda,\rho,\lambda',\rho',\iota)$$ and $$(G,e,i,m,\alpha,\lambda,\rho,\lambda',\rho',\iota'),$$ and that would be an
unintended departure from traditional mathematics.  If we were able to prove that $\iota = \iota'$, then our two groups would be equal.
Luckily, it is a theorem of Voevodsky that, for any type $X$, the type $isSet(X)$ is a proposition.\footnote{His proof depends on two axioms,
  both of which are consequences of the Univalence Axiom to be presented later.  One of them is a strong form of function extensionality, which
  states that functions having equal values are equal.}  Applied to $G$,
it shows that $\iota = \iota'$, implying in turn that our two groups are equal.  This fortuitous foundational result helps to show the feasibility of
the approach.

We now define $Group$ to be the type consisting of such 10-tuples; its elements are the groups living in $U_n$.  For precision about the
universe, one may use a subscript, as in $Group_n$.

It should now be evident to the reader how to follow the pattern established above and formalize the definitions of the types of many sorts of
mathematical objects, such as of Abelian groups, monoids, rings, commutative rings, and fields.  For example, the type $Set_n$ of all sets in
$U_n$ will consist of the pairs $(X,\iota)$, where $X$ is an element $U_n$, and $\iota$ is a proof that $X$ is a set.  As above, one sees that
two proofs that $X$ is a set yield equal elements of $Set_n$.

Now let's consider the problem of formalizing the notion of $G$-torsor.  By definition, it is a nonempty set $X$ with a free
transitive\footnote{One says that $G$ acts {\em freely} on $X$ if $gx=x$ implies $g=e$.  One says that $G$ acts {\em transitively} on $X$ if,
  for all $x,y \in X$, there is some $g \in G$ such that $gx=y$.  Thus $G$ acts freely and transitively on $X$ if, for every $x \in X$, the
  function of type $G \to X$ given by $g \mapsto g x$ is a bijection.  The conjunction of these two conditions is equivalent to the map $G
  \times X \to X \times X$ provided by $(g,x) \mapsto (gx,x)$ being a bijection, which is why these two conditions are usually stated together.}
left action by $G$.  For example, if $G$ is a subgroup of a group $H$, then the cosets of $G$ in $H$ are $G$-torsors.  We may represent a
$G$-torsor as a tuple $$(X,m,\alpha,\lambda,\iota,\tau,\nu),$$ where $m : G \times X \to X$ is the action, and the remaining components
represent the associativity property, the left unit property, the proof that $X$ is a set, the proof that $G$ acts freely and transitively on
$X$, and the proof that $X$ is nonempty\footnote{Nonemptiness of $X$ does not follow from the other properties.}.

Well, almost: what type does the proof $\nu$ of nonemptiness of $X$ belong to?  The only type available for the task seems to be $X$ itself:
suppose we were to use it.  Let $\nu$ and $\nu'$ be unequal elements of $X$.  Then we would get two unequal $G$-torsors,
$$(X,m,\alpha,\lambda,\iota,\tau,\nu)$$ and $$(X,m,\alpha,\lambda,\iota,\tau,\nu'),$$ and that would be a departure from traditional
mathematics: two proofs of nonemptiness of the set underlying a torsor should not give two different torsors.  Moreover, equipping a torsor with
a point renders the torsor canonically trivial.  In set theory there is a difference between a pointed set and a nonempty set, and we need a
formal way to express that difference in type theory.

In order to ensure that two proofs of nonemptiness (such as $\nu$ and $\nu'$) are always equal, we need them to be elements of a {\em
  proposition}, to be called $\nonempty{X}$, say, and to be called the {\em propositional truncation} of $X$.  Any element of $X$ should provide
a proof of nonemptiness, so there should be a map $\mu : X \to \nonempty{X}$.  One way to implement propositional truncation is to add it to the
language and justify it later, as is done in \citep{AwodeyBauer_prop_as_types}, or one could provide the following economical definition of
it, as is done in \citep{UniMath}.

For any type $X$ in $U_n$, Voevodsky defines $$\nonempty{X} \defeq \prod_{P:U_n} ( isProp(P) \to ((X \to P) \to P) ).$$
It is the conjunction of all the propositions ``implied'' by $X$, and thus
it is a proposition.  The map $\mu$ can be implemented as $x \mapsto P \mapsto i \mapsto f \mapsto f(x)$.  Moreover, it has a
universal property: any map $g : X \to Q$ to a proposition $Q$ in $U_n$ factors through $\mu$; the factorization map is provided by $ w \mapsto
((w(P))(j))(g)$, where $j$ is the proof that $Q$ is a proposition.\footnote{The definition is quantified over a variable $P$ in the universe
  $U_n$, so $\nonempty{X}$ lies not in $U_n$, as desired, but in $U_{n+1}$.
  Moreover, one wants to be able to apply the universal property to propositions $Q$ living in higher universes.
  To address these problems, Voevodsky's proposed system includes various ``resizing'' rules \citep{VV-resizing}, one of which ensures that
  propositions descend to lower universes.  He believes it likely that the resizing axioms are consistent with the rest of the system, but the
  simplicial set interpretation described below does not validate them.  Settling this issue is important.
}

Returning to our formalization above of the type of all $G$-torsors, we may now specify that the final component $\nu$ of the tuple is to be of
type $\nonempty{X}$, thereby completing the formalization.

Similarly, we may use propositional truncation to define the disjunction of two propositions as $P \vee Q \defeq \nonempty{P \amalg Q}$.  It is
a proposition and has the properties one would expect.

We may use propositional truncation to define the existential quantifier $$\exists _{x:X} P(x) \defeq \nonempty{ \sum _{x:X} P(x) }$$ for a family of
propositions $P(x)$ depending on a parameter $x:X$.  It expresses the statement that there is some element $x$ so that $P(x)$ is true, without
providing $x$.
It is a proposition and has the properties one would expect.

Surjectivity of a function $f : X \to Y$ may be encoded by the type $\forall _{y:Y} \nonempty{ f^{-1} (y) }$.  The type $\prod_{y:Y} f^{-1} (y)$
would not provide a suitable encoding, since it is the type of sections of $f$.

\section{Univalence}

At this point, we are faced with a design decision: whether to arrange for some of our identity types $a=b$ to have multiple elements.
Additional equalities would be useful, because they can be used in substitutions.  Alternatively, if classical mathematics were our only guide,
we would be led to introduce an axiom that asserts that every type is a set.  That wouldn't lead to a contradiction, but it wouldn't be
essential, because the classical constructions of new sets that don't involve Grothendieck universes, when formalized in type theory, yield
types that can be proved to be sets anyway, as Voevodsky has shown.

Let us consider the identity types $X=Y$ corresponding to types $X$ and $Y$ in the same universe, $U_n$ say.  The only obvious statement about
such types is that $X=X$ has an explicit and trivial proof: $\refl(X)$.

We motivate a desire for further proofs of equality between types as follows.

In set theory, the equation $\NN = \{ x \in \ZZ \mid x \ge 0 \}$ is false, if one uses the usual definition that $\ZZ$ is a certain quotient set
of $\NN \times \NN$.  The statement is false because the elements of the two sets are different.  The strongest thing that can be said is that
there is an isomorphism $\NN \weq \{ x \in \ZZ \mid x \ge 0 \}$.  Nevertheless, a mathematician may {\em identify} the two sets and
expect not to get into trouble.

In set theory, the associativity $(X \amalg Y) \amalg Z = X \amalg (Y \amalg Z)$ of binary sums is false, if one uses the usual definition that
$X \amalg Y$ is the set of pairs $(i,z)$ where $i = 0$ or $i=1$, $z \in X$ if $i = 0$, and $z \in Y$ if $i = 1$.  The statement is false because
the elements of the two sets are different.  The strongest thing that can be said is that there is a natural isomorphism
$(X \amalg Y) \amalg Z \weq X \amalg (Y \amalg Z)$.
Nevertheless, a mathematician may {\em identify} the two sets and expect not to get into trouble.

In type theory we may entertain the possibility that $\NN$ and $\{ x \in \ZZ \mid x \ge 0 \}$ are equal, and that binary sums are associative,
because the argument above, about the elements of one set not being equal to the elements of the other, is not available.  We may even try to
re-engineer our formal system to make those equalities happen.

Here's how Voevodsky made it happen.

A function $f : X \to Y$ between two types is called an {\em equivalence} if, for each $y:Y$, the fiber $f^{-1}(y)$ has exactly one element.  Voevodsky
has proven that being an equivalence is a proposition.\footnote{The proof depends on two axioms that are consequences of the Univalence
  Axiom below.}  The interested reader may encode this property as a type by applying encodings
previously introduced.  Thus the equivalences are the elements of a type, which is denoted by $X \weq Y$.  This definition has almost all of the
properties one would expect for the appropriate notion of isomorphism between two types.  For example, one can prove that the identity function
$X \to X$ is an equivalence, that an equivalence has an inverse function that is also an equivalence, and that a function with an inverse
function is an equivalence.

For types $X,Y:U_n$, we may define a function $ \Phi_{X,Y} : (X=Y) \to (X \weq Y)$ by induction: in the case where $Y$ is $X$, it sends
$\refl(X)$ to the identity equivalence.  This gives a relationship between equalities and equivalences.

Voevodsky's Univalence Axiom is stated as follows.

\begin{axiom}[Univalence]
  The map $ \Phi_{X,Y} : (X=Y) \to (X \weq Y)$ is an equivalence.
\end{axiom}

This axiom provides a way to promote equivalences into equalities, by application of the inverse of $ \Phi_{X,Y} $, thereby satisfying the
desires motivated above.\footnote{The statement of the Univalence Axiom is (encoded by a type that is) a proposition, because being an equivalence
  is a proposition.  Mart\'in Escard\'o has pointed out that the proof can easily be made independent of any axioms, because, while proving
  equality of any two elements of the type asserting univalence, one has available all the consequences of univalence.  It's desirable for
  statements of axioms to be propositions, so asserting the axioms doesn't lead to indeterminate behavior.}

Now we address an issue of terminology.  From this point on, the reader may prefer to think of an element of $X=Y$ not as providing a proof that
the types $X$ and $Y$ are {\em equal}, but as providing a way to {\em identify} $X$ with $Y$.  Different elements will provide different ways.
This alternative reading may relieve some discomfort for readers who prefer to preserve the word ``equality'' for something that can happen in
just one way.  In line with that consideration, we will refer to the type $X=Y$ as an {\em identity type}, rather than as an {\em equality
  type}, and we will refer to its elements as {\em identities} or {\em identifications}.
\par
For example, a bijection $\theta : \NN \weq \{ x \in \ZZ \mid x \ge 0 \}$ provides, via univalence, a proof $\theta'$ that
$\NN = \{ x \in \ZZ \mid x \ge 0 \}$, which in turn allows us to use $\theta'$ to {\em identify} $\NN$ with $\{ x \in \ZZ \mid x \ge 0 \}$.
This way of putting it is common mathematical practice.  What's new is support for the practice in the language of logic.
\par
Let us pause to infer an important foundational consequence.  We know that the statements of type theory are invariant under identity (by
substitution).  With the ability to promote equivalences to identities, we see that the statements of type theory are invariant under
equivalence, thereby fulfilling the mathematicians' dream: one cannot express a property in our formal language that fails to be invariant under
equivalence.  That invariance is a strong and useful principle, which flows from the precise way that the formal language has been restricted to
the minimum required for formalization of mathematics.

Pairs of sets often have more than one isomorphism between them.  Hence, in the presence of the univalence axiom, pairs of types often will have
more than one proof of identity, and thus the universes containing those types are not sets.  Moreover, in the absence of univalence and of
axioms contradicting it, the
universes cannot be proved to be sets.  This is a major difference from set theory.

A consequence of the Univalence Axiom is that isomorphism between groups is the same thing as identity.  In other words, given two groups $G, H
: Group$, let $G \isom H$ denote the type of isomorphisms between $G$ and $H$.  Then the natural map $G = H \to G \isom H$ is an equivalence.  An
analogous statement can be proved for all the other standard mathematical structures: partially ordered set, $G$-torsor, monoid, ring, etc.
This general principle is called the {\em structure identity principle} in the book \citep[section 9.8]{hottbook}.  The main point in the proof
of it is that the notion of isomorphism between two instances of an algebraic structure pays attention to every element of the structure, so the
elements of either structure are uniquely determined by the equivalence on the underlying types and the elements of the other structure.  As a
consequence, we see that $Group$ is not a set.

The structure identity principle can be used to prove, for example, that if $X$ and $Y$ are partially ordered sets, $x$ is a minimal element of
$X$, and $f : X \weq Y$ is an isomorphism of partially ordered sets, then $f(x)$ is a minimal element of $Y$.  The same proof works for maximal
elements: neither proof needs to appeal to the definition of ``minimal element'' or ``maximal element'', because the principle applies to all
properties expressible in our formal language.  By contrast, in set theory, to prove those assertions, the definitions of minimal element and of
maximal element must be examined and used, since the more permissive language of set theory permits statements to be written that contain
mathematical irrelevancies; for an example, consider the proposition $x=2$, which does not imply $f(x)=2$.

Traditionally, the objects of a category have constituted a set, but the natural way to formalize the {\em category} of groups is with its
objects being the elements of the {\em type} $Group$.  The best categories are the ones (such as this one) where the identities in the type of
objects are equivalent to the isomorphisms, for it is such categories that fulfill the category theorists' dream: to work in a mathematical
language where one cannot express a property that fails to be invariant under isomorphism\footnote{For example, Makkai \citep{MR1678360} has
  formulated the dream this way, as a requirement to be satisfied by a future ``Structuralist Foundation of Abstract Mathematics'' (SFAM):
  ``SFAM adopts isomorphism of objects in a category ... to play the role of equality for those objects, in the sense that ... all grammatically
  correct properties of objects of a fixed category are to be invariant under isomorphism.''}; such categories are called {\em univalent}.  See
\citep{MR3340533} for a procedure that renders categories univalent in a universal way.
\par
On the other hand, sometimes one needs a category whose objects form a set; for example, one may wish
to construct the geometric realization of the category as a topological space.  If so, then one may incorporate enough extraneous data in each object of the
category to make proofs of identity between two objects unique.  For example, one may equip each object in the category of finite dimensional
vector spaces with an ordered basis -- that is enough, because isomorphisms respecting the bases are unique.

As we said above, in the presence of univalence the type $Group$ is not a set, but there is something positive to be said about it.  Given two
groups $G, H : Group$, the type $G = H$ is a set, because the map $G = H \to G \isom H$ is an equivalence of it with a set.
\par
This motivates another fundamental definition of Voevodsky's, which unifies and simplifies the formalization of the basic theorems of the
subject.  We define by induction on a natural number $n$ what it means to say that a type $X$ has {\em $h$-level (at most) $n$}: (1) it has
$h$-level $0$ if it has exactly one element; (2) it has $h$-level $n+1$ if for any elements $a$ and $b$ of $X$, the type $a=b$ has $h$-level
$n$.  Voevodsky has proven that the types of $h$-level 1 are the propositions\footnote{Actually, he {\em defined} a proposition to be a type of
  $h$-level $1$; our definition differs slightly.}, and then one can see from the
definition that the types of 
$h$-level 2 are the sets.  We see also that $Group$ is of $h$-level $3$, which is the natural level for the type of objects of any groupoid.
Properly defined, the type of categories will be of $h$-level 4, because the groupoid of categories is a 2-groupoid: an equivalence between two
categories may have automorphisms.  In general one expects the objects of an {\em $n$-groupoid} to be of $h$-level $n+2$.

Voevodsky points out a stratification of all of mathematics into levels that follows from the stratification of types by $h$-level: {\em
  element-level} mathematics concerns equations between elements of sets, or properties of elements of sets; {\em set-level} mathematics
concerns isomorphisms between algebraic structures, or properties of algebraic structures invariant under isomorphism; {\em groupoid-level}
mathematics concerns equivalences of groupoids, or properties of groupoids invariant under equivalence of groupoids; and so on
\citep{UniMath2015}.

Another common mathematical practice is to speak of ``natural'' or ``canonical'' constructions; for example, one often states that there is no
natural isomorphism between a vector space $V$ and its dual vector space $V^*$.  This practice is also directly supported now by logic, provided
we consider naturality only with respect to isomorphisms $V \isom V'$, which, by the structure identity principle, are captured by identities,
which in turn may be used in substitutions.  We may formalize the type of isomorphisms between a vector space and its dual as the type
$\prod_{V:Vect} ( V \isom V^* )$, where $Vect$ is the type of vector spaces over a field $k$ in universe $U_n$.  In the presence of
univalence, one may prove that the type of such families of isomorphisms is empty.\footnote{Compare with the discussion in \citep{mcall} or the
  use of the word ``natural'' in the proof of \citep[Theorem 3.2.2]{hottbook}.}

\section{The interpretation}
\label{the_interpretation}

The possibility that an identity type $a=b$ has multiple elements presents a problem for the mathematical reader who wants to understand the
situation, as well as a problem for the reader who has perused the interpretation for our formal system in classical mathematics that we
sketched above.  That interpretation is incompatible with the Univalence Axiom, because, roughly speaking, for a set $X:U_n$ with more than one
element it interprets the map $ \Phi_{X,X} : (X=X) \to (X \weq X)$ as a function from the one point set to the set of permutations of $|X|$; the
function sends the single point $*$ to the identity bijection.  The function is not a bijection, as would be required to interpret the
Univalence Axiom.

The idea for repairing this is to interpret each type $X$ as a topological space\footnote{The topological spaces used in the interpretation are
  actually fibrant simplicial sets.} $|X|$, to interpret an element $a:X$ as a point $|a|$ of the space $|X|$, to interpret a proof of $a=b$ as
a path connecting $|a|$ to $|b|$ in $|X|$, and to interpret the type $a=b$ as the space of all such paths.  Symmetry of identity passes to
reversal of paths in the interpretation, and transitivity of identity passes to concatenation of paths.  A family $Y(x)$ of types parametrized
by a variable $x:X$ will be interpreted by a fibration\footnote{In topology, a ``continuous'' family of topological spaces is one that is provided
  by the fibers of a continuous map that is a fibration.} $E(Y)\to |X|$ whose fibers over the points $|a| \in |X|$ arising from elements
$a:X$ are the spaces $|Y(a)|$.  Induction for identity will be interpreted by the lifting property that characterizes fibrations, as observed in
\citep{0709.0248}.  One assumes one is given an ascending sequence of Grothendieck universes $V_0 \in V_1 \in \dots$ and one interprets the
universe $U_n$ by the space $|U_n|$ of all spaces in $V_n$; it is a space in $V_{n+1}$.  For $n \ge 0$, types of $h$-level $n+2$ are interpreted
by spaces whose homotopy groups vanish above dimension $n$.  Validity of the Univalence Axiom in this interpretation is a theorem of Voevodsky \citep{1203.2553},
which builds on a 2006 construction of A.{} Bousfield to show that the path lifting map from the space of paths between two
points of $|U_n|$ to the space of homotopy equivalences between the
corresponding fibers of the universal space over $|U_n|$ is a homotopy equivalence.  The expository paper \citep{1211.2851} gives the main ideas
of the construction of the interpretation; full details of the correctness of the interpretation, and thus of the soundness of the theory, will appear in a
series of papers planned by Voevodsky, aimed at treating a large class of formal languages, some of which are available: \citep{39, 42, 40, 104,
  103, 109, 110, 112, MR3607210, MR3607209}.  What necessitates so many details to be expressed in print is the abundance of grammatical rules
in a formal language such as this one, each one of which has an incarnation in the interpretation that must be carefully considered, along with
the way they interact with each other.

The Law of Excluded Middle and the Axiom of Choice are also validated by this interpretation, so classical mathematics is supported soundly by
the univalent foundations.

\section{Further developments}\label{otherdev}

In addition to the focus of many on establishing the soundness of the system and on paving the way for widespread adoption of univalent
foundations and the use of proof assistants by mathematicians, there are other interesting developments focusing such things as computability
and the study of types in their own right.

The phrase ``synthetic homotopy theory'' refers to the enterprise where one regards a type as a good substitute for the classical notion of
topological space, and one regards a proof of identity as a good substitute for the notion of a path between two points of a space: the goal is
to see which theorems of homotopy theory have analogues that remain provable in this context.  This context is a primitive one, because so few
axioms are assumed to set up the theory, as we have seen, so the theorems that hold in this context will the most fundamental ones, capable of
the most generalization.  An amazing surprise is that so many theorems of homotopy theory hold in this rarefied context.

Given a type $X$ we define the type $\pi_0 X$ of connected components as the quotient of $X$ by the equivalence relation $b=c$, where $b$ and
$c$ are elements of $X$.  Voevodsky has proven that $\pi_0 X$ is a set (as is the case for the type of equivalence classes of any equivalence
relation on $X$).

Given a type $X$ and a basepoint $a:X$, we define the loop space $\Omega X$ to be the type $a=a$, and we equip it with the basepoint $\refl(a) : \Omega X$.
Iterating $n$ times yields a type $\Omega^n X$, and we may define $\pi_n X := \pi_0 \Omega^n X$.  The proofs of transitivity and symmetry in
section \ref{paths} provide $\pi_1 X$ with a group structure, and thus provide $\pi_n X$ with $n$ group structures when $n > 0$.  A basic fact
of homotopy theory is that the group $\pi_n X$ is abelian for $n \ge 2$, and that the various group structures are the same.  To prove that, it
suffices to show that the two composition operations on the type $\Omega^2 X$ (coming from transitivity) agree and are commutative.  (We don't
have to say ``commutative up to homotopy'' here, because in our context, homotopy and identity are the same.)  One standard argument, due to
Eckmann and Hilton, involves showing that a monoid object in the category of monoids is commutative and the two operations coincide.  That
argument is formal enough that it applies here \citep[Theorem 2.1.6]{hottbook} and gives the first indication that something wonderful is going
on.  The interested reader may refer to \citep[Chapter 8]{hottbook} or to \citep{1703.03007} for more leisurely expositions of this pursuit.

It is not yet known how to construct types that correspond to the spheres $S^n$ in the univalent foundations as formalized by Voevodsky (except for $S^1$, which can be
defined as the type $B\ZZ$ of $\ZZ$-torsors).  The book \citep{hottbook} introduces spheres by adding further basic type constructors to the
formal language -- they are called ``higher inductive types''.  For example, the circle $S^1$ is defined inductively by declaring that there is
a basepoint $s : S^1$ and that there is a loop $\ell : s = s$.  The higher spheres are done in a similar way.  Voevodsky's program for proving
soundness of univalent foundations doesn't include higher inductive types, but others are pursuing it: see \citep{1705.07088}.  Much of the work
cited below uses higher inductive types.

In \citep{MR3323808} one may find a proof that the fundamental group of $S^1$ is $\ZZ$, and in \citep{MR3566711} one may find a formalization of
the Seifert--van Kampen theorem.
In \citep{1605.03227} one may find a formalization of the Blakers-Massey connectedness theorem, using higher inductive types to construct
homotopy pushouts.  A translation of the formal proof into classical homotopy theory was given in \citep{RezkBM}, with the expectation that it
would go through for any $\infty$-topos.  Finally, the proof was (manually) translated into the language of $\infty$-topoi in the paper \citep{1703.09050},
yielding a new theorem.\footnote{The companion paper \citep{1703.09632} translates the proof into the language of Goodwillie's calculus of
  functors (using a different ``modality'') and yields a new result in that context, which has various known results as corollaries.}
Even better would be a meta-theorem that says that type theory serves as an internal language for higher topoi, so manually rewriting and
rechecking the proof in the new context would not be required.  Some progress on that dream has recently been made in the papers \citep{1610.00037,1709.09519}.

A type-theoretic proof that $\pi_4 (S^3) \cong \ZZ/2$ is offered by the thesis \citep{1606.05916} and the paper \citep{1710.10307}.
The first part of the proof demonstrates the
existence of a natural number $n$ satisfying $\pi_4 (S^3) \cong \ZZ/n$, and the second part demonstrates that $n=2$.  The proofs are
constructive, in the sense that the use of axioms (such as the Law of Excluded Middle or the Axiom of Choice) is avoided, aside from the use of
the Univalence Axiom.  Although normally the use of any axiom interferes with computability, Voevodsky conjectured (albeit for a system not
including higher inductive types) that the interference arising from the use of the Univalence Axiom can be bypassed.  Thus one may hope to
design a proof assistant that can produce the value of $n$ by performing a computation.  Fundamental progress has been made in that direction
through the development of ``cubical type theory'' \citep{MR3281415,1611.02108,1802.01170} and of a proof assistant, {\tt cubicaltt} \citep{cubicaltt}, based
on it.  See also the continuing development of {\em RedPRL} \citep{RedPRL}, described in the papers
\citep{harper-et-al-1,harper-et-al-2,harper-et-al-3,harper-et-al-4}, and see also \citep{cart-cube}.
A formal expression for the number $n$ has been submitted to {\tt cubicaltt} for
evaluation, but the computation ran out of memory after running for 6 hours on a machine with plenty of memory, according to Brunerie.  Work in
that direction continues.

In \citep{1509.07584} one may find a type theoretic proof of the Brouwer fixed point theorem, accomplished by an enhanced type theory called
{\em real-cohesive homotopy type theory}, which supports synthetic topology side by side with synthetic homotopy theory.
In \citep{1802.02191} is a synthetic proof of the theorem that cellular cohomology computes ordinary cohomology of CW-complexes, under the
assumption that the CW-complexes are finite.

\bibliographystyle{amsplain}
\bibliography{papers}

\providecommand{\bysame}{\leavevmode\hbox to3em{\hrulefill}\thinspace}
\providecommand{\MR}{\relax\ifhmode\unskip\space\fi MR }
\providecommand{\MRhref}[2]{%
  \href{http://www.ams.org/mathscinet-getitem?mr=#1}{#2}
}
\providecommand{\href}[2]{#2}
\begin{thebibliography}{10}

\bibitem{HoTT-Agda}
\emph{{Homotopy Type Theory in Agda}}, an Agda library of formalized proofs,
  available at \url{https://github.com/HoTT/HoTT-Agda}.

\bibitem{Lean}
\emph{Lean theorem prover}, a proof assistant together with a library of
  formalized proofs, available at \url{https://github.com/leanprover/lean}.

\bibitem{RedPRL}
\emph{{RedPRL}}, a proof assistant for Computational Cubical Type Theory
  together with a library of formalized proofs, available at
  \url{http://www.redprl.org/}.

\bibitem{HoTT}
\emph{{The HoTT Library}}, a Coq library of formalized proofs, available at
  \url{https://github.com/HoTT/HoTT}.

\bibitem{MR3340533}
Benedikt Ahrens, Krzysztof Kapulkin, and Michael Shulman, \emph{Univalent
  categories and the {R}ezk completion}, Math. Structures Comput. Sci.
  \textbf{25} (2015), no.~5, 1010--1039. \MR{3340533}

\bibitem{1703.09050}
Mathieu Anel, Georg Biedermann, Eric Finster, and Andr\'e Joyal, \emph{{A
  Generalized Blakers-Massey Theorem}}, 2017, \arxiv{1703.09050}.

\bibitem{1703.09632}
\bysame, \emph{{Goodwillie's Calculus of Functors and Higher Topos Theory}},
  2017, \arxiv{1703.09632}.

\bibitem{cart-cube}
Carlo Angiuli, Guillaume Brunerie, Thierry Coquand, Kuen-Bang Hou, Robert
  Harper, and Daniel~R. Licata, \emph{{Cartesian Cubical Type Theory}}, 2017,
  {a}vailable at
  \url{https://github.com/dlicata335/cart-cube/raw/master/cart-cube.pdf}.

\bibitem{harper-et-al-2}
Carlo Angiuli and Robert Harper, \emph{{Computational higher type theory II:
  Dependent cubical realizability}}, 2016, \arxiv{1606.09638}.

\bibitem{harper-et-al-1}
Carlo Angiuli, Robert Harper, and Todd Wilson, \emph{{Computational higher type
  theory I: Abstract cubical realizability}}, 2016, \arxiv{1604.08873}.

\bibitem{harper-et-al-3}
Carlo Angiuli, Kuen-Bang Hou, and Robert Harper, \emph{{Computational Higher
  Type Theory III: Univalent Universes and Exact Equality}}, 2017,
  \arxiv{1712.01800}.

\bibitem{MR0543792}
K.~Appel and W.~Haken, \emph{Every planar map is four colorable. {I}.
  {D}ischarging}, Illinois J. Math. \textbf{21} (1977), no.~3, 429--490.
  \MR{0543792}

\bibitem{MR0543793}
K.~Appel, W.~Haken, and J.~Koch, \emph{Every planar map is four colorable.
  {II}. {R}educibility}, Illinois J. Math. \textbf{21} (1977), no.~3, 491--567.
  \MR{0543793}

\bibitem{10.2307/986491}
Kenneth~I. Appel, \emph{The use of the computer in the proof of the four color
  theorem}, Proceedings of the American Philosophical Society \textbf{128}
  (1984), no.~1, 35--39.

\bibitem{0709.0248}
Steve Awodey and Michael~A. Warren, \emph{Homotopy theoretic models of identity
  types}, Math. Proc. Cambridge Philos. Soc. \textbf{146} (2009), no.~1,
  45--55, \arxiv{0709.0248}. \MR{2461866}

\bibitem{AwodeyBauer_prop_as_types}
Steven Awodey and Andrej Bauer, \emph{Propositions as [types]}, J. Logic
  Comput. \textbf{14} (2004), no.~4, 447--471. \MR{2081047}

\bibitem{1610.04591}
Andrej Bauer, Jason Gross, Peter~LeFanu Lumsdaine, Mike Shulman, Matthieu
  Sozeau, and Bas Spitters, \emph{{The HoTT Library: A formalization of
  homotopy type theory in Coq}}, 2016, \arxiv{1610.04591}.

\bibitem{IMGelfand}
A.~Beilinson, \emph{I. {M}. {G}elfand and his seminar---a presence}, Notices
  Amer. Math. Soc. \textbf{63} (2016), no.~3, 295--298, \arxiv{1505.00710}.
  \MR{3445168}

\bibitem{bernays-1}
Paul Bernays, \emph{A system of axiomatic set theory. {P}art {I}}, J. Symbolic
  Logic \textbf{2} (1937), 65--77.

\bibitem{MR3281415}
Marc Bezem, Thierry Coquand, and Simon Huber, \emph{A model of type theory in
  cubical sets}, 19th {I}nternational {C}onference on {T}ypes for {P}roofs and
  {P}rograms, LIPIcs. Leibniz Int. Proc. Inform., vol.~26, Schloss Dagstuhl.
  Leibniz-Zent. Inform., Wadern, 2014, pp.~107--128. \MR{3281415}

\bibitem{1606.05916}
Guillaume Brunerie, \emph{{On the homotopy groups of spheres in homotopy type
  theory}}, 2016, {Ph.D.} thesis, \arxiv{1606.05916}.

\bibitem{1710.10307}
\bysame, \emph{{The James construction and $\pi_4(\mathbb{S}^3)$ in homotopy
  type theory}}, 2017, \arxiv{1710.10307}, accepted for publication in the
  Journal of Automated Reasoning.

\bibitem{1802.02191}
Ulrik Buchholtz and Kuen-Bang~(Favonia) Hou, \emph{{Cellular Cohomology in
  Homotopy Type Theory}}, February, 2018, \arxiv{1802.02191}.

\bibitem{harper-et-al-4}
Evan Cavallo and Robert Harper, \emph{{Computational Higher Type Theory IV:
  Inductive Types}}, 2018, \arxiv{1801.01568}.

\bibitem{cubicaltt}
Cyril Cohen, Thierry Coquand, Simon Huber, and Anders M\"ortberg,
  \emph{Experimental implementation of cubical type theory}, the proof
  assistant {\tt cubical}, written in Haskell, available at
  \url{https://github.com/mortberg/cubicaltt}.

\bibitem{1611.02108}
\bysame, \emph{{Cubical Type Theory: a constructive interpretation of the
  univalence axiom}},  (2015), To appear in the post-proceedings of TYPES 2015.
  Preprint avaiable at https://arxiv.org/abs/1611.02108.

\bibitem{MR3363596}
Thierry Coquand, \emph{Th\'eorie des types d\'ependants et axiome
  d'univalence}, Ast\'erisque (2015), no.~367-368, Exp. No. 1085, x, 367--386.
  \MR{3363596}

\bibitem{MR0274219}
N.~G. de~Bruijn, \emph{The mathematical language {AUTOMATH}, its usage, and
  some of its extensions}, Symposium on {A}utomatic {D}emonstration
  ({V}ersailles, 1968), Lecture Notes in Mathematics, Vol. 125. Springer,
  Berlin, 1970, pp.~29--61. \MR{0274219}

\bibitem{MR2463991}
Georges Gonthier, \emph{Formal proof---the four-color theorem}, Notices Amer.
  Math. Soc. \textbf{55} (2008), no.~11, 1382--1393. \MR{2463991}

\bibitem{10.1007/978-3-642-39634-2_14}
Georges Gonthier, Andrea Asperti, Jeremy Avigad, Yves Bertot, Cyril Cohen,
  Fran{\c{c}}ois Garillot, St{\'e}phane Le~Roux, Assia Mahboubi, Russell
  O'Connor, Sidi Ould~Biha, Ioana Pasca, Laurence Rideau, Alexey Solovyev,
  Enrico Tassi, and Laurent Th{\'e}ry, \emph{A machine-checked proof of the odd
  order theorem}, Interactive Theorem Proving (Berlin, Heidelberg) (Sandrine
  Blazy, Christine Paulin-Mohring, and David Pichardie, eds.), Springer Berlin
  Heidelberg, 2013, pp.~163--179.

\bibitem{kepler}
Thomas Hales, Mark Adams, Gertrud Bauer, Tat~Dat Dang, John Harrison, Le~Truong
  Hoang, Cezary Kaliszyk, Victor Magron, Sean Mclaughlin, Tat Thang Nguyen,
  et~al., \emph{{A formal proof of the Kepler conjecture}}, Forum of
  Mathematics, Pi \textbf{5} (2017), e2.

\bibitem{MR3566711}
Kuen-Bang Hou and Michael Shulman, \emph{The {S}eifert--van {K}ampen theorem in
  homotopy type theory}, Computer science logic 2016, LIPIcs. Leibniz Int.
  Proc. Inform., vol.~62, Schloss Dagstuhl. Leibniz-Zent. Inform., Wadern,
  2016, pp.~Art. No. 22, 16. \MR{3566711}

\bibitem{1605.03227}
Kuen-Bang Hou~(Favonia) Hou, Eric Finster, Dan Licata, and Peter~LeFanu
  Lumsdaine, \emph{{A mechanization of the Blakers-Massey connectivity theorem
  in Homotopy Type Theory}}, 2016, \arxiv{1605.03227}.

\bibitem{1211.2851}
Chris Kapulkin and Peter~LeFanu Lumsdaine, \emph{The {S}implicial {M}odel of
  {U}nivalent {F}oundations (after {V}oevodsky)}, 2012, \arxiv{1211.2851}.

\bibitem{1610.00037}
\bysame, \emph{{The homotopy theory of type theories}}, 2016,
  \arxiv{1610.00037}.

\bibitem{1203.2553}
Chris Kapulkin, Peter~LeFanu Lumsdaine, and Vladimir Voevodsky,
  \emph{Univalence in simplicial sets}, 2012, \arxiv{1203.2553}.

\bibitem{1709.09519}
Chris Kapulkin and Karol Szumiło, \emph{{Internal Language of Finitely
  Complete $(\infty,1)$-categories}}, 2017, \arxiv{1709.09519}.

\bibitem{MR3323808}
Daniel~R. Licata and Michael Shulman, \emph{Calculating the fundamental group
  of the circle in homotopy type theory}, 2013 28th {A}nnual {ACM}/{IEEE}
  {S}ymposium on {L}ogic in {C}omputer {S}cience ({LICS} 2013), IEEE Computer
  Soc., Los Alamitos, CA, 2013, pp.~223--232. \MR{3323808}

\bibitem{1705.07088}
Peter~LeFanu Lumsdaine and Mike Shulman, \emph{{Semantics of higher inductive
  types}}, 2017, \arxiv{1705.07088}.

\bibitem{MR1678360}
M.~Makkai, \emph{Towards a categorical foundation of mathematics}, Logic
  {C}olloquium '95 ({H}aifa), Lecture Notes Logic, vol.~11, Springer, Berlin,
  1998, pp.~153--190. \MR{1678360}

\bibitem{MR0387023}
Per Martin-L{\"o}f, \emph{Hauptsatz for the intuitionistic theory of iterated
  inductive definitions}, Proceedings of the {S}econd {S}candinavian {L}ogic
  {S}ymposium ({U}niv. {O}slo, {O}slo, 1970), North-Holland, Amsterdam, 1971,
  pp.~179--216. Studies in Logic and the Foundations of Mathematics, Vol. 63.
  \MR{0387023 (52 \#7870)}

\bibitem{MR0387009}
\bysame, \emph{An intuitionistic theory of types: predicative part}, Logic
  {C}olloquium '73 ({B}ristol, 1973), North-Holland, Amsterdam, 1975,
  pp.~73--118. Studies in Logic and the Foundations of Mathematics, Vol. 80.
  \MR{0387009 (52 \#7856)}

\bibitem{MLTT79}
Per Martin-L\"of, \emph{Constructive mathematics and computer programming},
  Logic, methodology and philosophy of science, {VI} ({H}annover, 1979), Stud.
  Logic Found. Math., vol. 104, North-Holland, Amsterdam, 1982, pp.~153--175.
  \MR{682410}

\bibitem{mcall}
David McAllester, \emph{Morphoid type theory, a typed platonic foundation for
  mathematics}, \arxiv{1407.7274}.

\bibitem{peano-principia}
Ioseph Peano, \emph{Arithmetices principia}, Fratres Bocca, 1889.

\bibitem{RezkBM}
Charles Rezk, \emph{{Proof of the Blakers-Massey Theorem}}, Available at
  \url{https://faculty.math.illinois.edu/~rezk/freudenthal-and-blakers-massey.pdf}.

\bibitem{MR1506041}
Bertrand Russell, \emph{Mathematical {L}ogic as {B}ased on the {T}heory of
  {T}ypes}, Amer. J. Math. \textbf{30} (1908), no.~3, 222--262. \MR{1506041}

\bibitem{Swansea-Shulman}
Michael Shulman, \emph{{Minicourse on Homotopy Type Theory}}, {slides from a
  talk, available} at
  \url{http://home.sandiego.edu/~shulman/hottminicourse2012/}, 2012.

\bibitem{1509.07584}
Michael Shulman, \emph{{Brouwer's fixed-point theorem in real-cohesive homotopy
  type theory}}, 2015, \arxiv{1509.07584}.

\bibitem{1703.03007}
Michael Shulman, \emph{{Homotopy type theory: the logic of space}}, 2017,
  \arxiv{1703.03007}.

\bibitem{Coq}
The Coq~Development Team, \emph{The coq proof assistant}, Available at
  \url{https://coq.inria.fr/}.

\bibitem{1802.01170}
Coquand Thierry, Simon Huber, and Anders M\"ortberg, \emph{{On Higher Inductive
  Types in Cubical Type Theory}}, 2018, \arxiv{1802.01170}.

\bibitem{Tsem4}
D.~Tsementzis, \emph{A meaning explanation for hott},  (2016), Available at
  \url{http://philsci-archive.pitt.edu/12824/}.

\bibitem{Tsem1}
\bysame, \emph{Univalent foundations as structuralist foundations}, Synthese
  \textbf{194} (2017), no.~9, 3583--3617, Available at
  \url{https://link.springer.com/article/10.1007/s11229-016-1109-x}.

\bibitem{Tsem2}
\bysame, \emph{What is a higher-level set?†}, Philosophia Mathematica (2017),
  (forthcoming), {Available at
  \url{https://academic.oup.com/philmat/article-abstract/2895117/What-is-a-Higher-Level-Set}}.

\bibitem{Tsem3}
D.~Tsementzis and H.~Halvorson, \emph{Foundations and philosophy},
  Philosopher's Imprint (forthcoming) (2017), Available at
  \url{http://philsci-archive.pitt.edu/13504/}.

\bibitem{hottbook}
The {Univalent Foundations Program}, \emph{Homotopy type theory: Univalent
  foundations of mathematics}, \url{https://homotopytypetheory.org/book},
  Institute for Advanced Study, 2013.

\bibitem{Foundations}
Vladimir Voevodsky, \emph{{Foundations: Development of the univalent
  foundations of mathematics in Coq}}, a Coq library of formalized proofs,
  available at \url{https://github.com/vladimirias/Foundations}.

\bibitem{VV-resizing}
\bysame, \emph{{Resizing Rules - their use and semantic justification}},
  {s}lides from a talk available at
  \url{https://www.math.ias.edu/vladimir/sites/math.ias.edu.vladimir/files/2011_Bergen.pdf},
  2011.

\bibitem{40}
\bysame, \emph{B-systems}, 2014, \arxiv{1410.5389}.

\bibitem{42}
\bysame, \emph{C-system of a module over a monad on sets}, 2014,
  \arxiv{1407.3394}.

\bibitem{VV-IHP-talk}
\bysame, \emph{Univalent foundations - new type-theoretic foundations of
  mathematics}, {s}lides from a talk at IHP, Paris, on April 22, 2014,
  available at
  \url{https://www.math.ias.edu/vladimir/sites/math.ias.edu.vladimir/files/2014_04_22_slides.pdf},
  2014.

\bibitem{109}
\bysame, \emph{A {C}-system defined by a universe category}, Theory Appl.
  Categ. \textbf{30} (2015), no.~37, 1181--1215, \arxiv{1409.7925}.
  \MR{3402489}

\bibitem{UniMath2015}
\bysame, \emph{An experimental library of formalized mathematics based on the
  univalent foundations}, Math. Structures Comput. Sci. \textbf{25} (2015),
  no.~5, 1278--1294. \MR{3340542}

\bibitem{103}
\bysame, \emph{Martin-{L}\"of identity types in the {C}-systems defined by a
  universe category}, 2015, \arxiv{1505.06446}, pp.~1--51.

\bibitem{104}
\bysame, \emph{Products of families of types in the {C}-systems defined by a
  universe category}, 2015, \arxiv{1503.07072}.

\bibitem{110}
\bysame, \emph{Lawvere theories and {J}f-relative monads}, arXiv:1601.02158
  (2016), 1--21, \arxiv{1601.02158}.

\bibitem{112}
\bysame, \emph{Products of families of types and {$(\Pi,\lambda)$}-structures
  on {C}-systems}, Theory Appl. Categ. \textbf{31} (2016), Paper No. 36,
  1044--1094. \MR{3584698}

\bibitem{39}
\bysame, \emph{Subsystems and regular quotients of {C}-systems}, A panorama of
  mathematics: pure and applied; Conference on Mathematics and its
  Applications, (Kuwait City, 2014), Contemp. Math., vol. 658, Amer. Math.
  Soc., Providence, RI, 2016, pp.~127--137. \MR{3475277}

\bibitem{MR3607209}
\bysame, \emph{C-systems defined by universe categories: presheaves}, Theory
  Appl. Categ. \textbf{32} (2017), Paper No. 3, 53--112. \MR{3607209}

\bibitem{MR3607210}
\bysame, \emph{The {$(\Pi, \lambda)$}-structures on the {C}-systems defined by
  universe categories}, Theory Appl. Categ. \textbf{32} (2017), Paper No. 4,
  113--121. \MR{3607210}

\bibitem{UniMath}
Vladimir Voevodsky, Benedikt Ahrens, Daniel Grayson, et~al., \emph{{UniMath ---
  a computer-checked library of univalent mathematics}}, {available} at
  \url{http://UniMath.org}.

\end{thebibliography}

\end{document}